\documentclass[11pt]{article}
\usepackage{amsmath}
\newcommand{\cn}{\mathrm{cn}}
\newcommand{\nc}{\mathrm{nc}}
\newcommand{\sd}{\mathrm{sd}}
\newcommand{\ds}{\mathrm{ds}}
\newtheorem{example}{\normalfont\em{Example}}

\title{Exponential--Weierstrass type, exponential--Jacobi type and solitary type
 solutions to some conformable fractional equations}

\author{Sirendaoreji\\College of Mathematical Science,Inner Mongolia Normal
University,\\ Huhhot 010022,Inner Mongolia, P.R. China.}
\begin{document}
\date{}
\maketitle
\begin{abstract}
A new algebraic method to find two special types of exact
traveling wave solutions and the solitary type solutions
to some conformable fractional partial differential equations
is proposed.
The two special types of solutions given by the product of
exponential and Weierstrass elliptic functions, and the product of
exponential and Jacobi elliptic functions are new and they cannot
be obtained by using the existing algebraic methods.\par
 \end{abstract}

\section{Introduction}\label{sec:1}

In recent years,the problem of finding exact solutions to nonlinear
fractional partial differential equations (FPDEs) has become an
active research area in nonlinear science due to the wide applicability
of nonlinear FPDEs to describe many important phenomena and dynamic processes
in fluid mechanics,biological and chemical processes, signal processing,
control systems and so on.
Based on the product rule and the chain rule for fractional derivatives,
several direct methods have been applied to find the exact solutions to
some nonlinear FPDEs\cite{Abdel, Ali,Feng, Gepreel,Guner}. However,the
recent studies showed that the fractional derivatives, such
as the Riemann--Liouville derivative and the Caputo derivative etc.
do not satisfy the product rule and the chain rule\cite{LiuCS15, LiuCS18}.
Therefore,the previous works of using these fractional derivatives
to solve those nonlinear FPDEs were not correct.
But a so-called $\alpha$--order conformable fractional derivative (CFD)
of a function $f: \left[0,\infty\right)\rightarrow{{R}}$ given by
\begin{equation}
\label{eq1:1}
D_t{^\alpha}f(t)=\lim_{\varepsilon\rightarrow{0}}{\frac {f(t+\varepsilon{t^{1-\alpha}})-f(t)}{\varepsilon}},
t>0,\alpha\in\left(0,1\right],
\end{equation}
which was defined by Khalil {\it et al}.~\cite{Khalil} just satisfies
the required properties such as~\cite{Khalil,Neirameh}
\begin{eqnarray}
\label{eq1:2}
&D_t^\alpha(af+bg)=aD_t^\alpha{f}+bD_t^\alpha{g},\forall{x}\in{R},
\mbox{(Linearity)}~\\
\label{eq1:3}
&D_t^\alpha{t^\beta}=\beta{t^{\beta-\alpha}},\\
\label{eq1:4}
&D_t^\alpha(fg)=gD_t^\alpha{f}+fD_t^\alpha{g},\mbox{(Product rule)},\\
\label{eq1:5}
&D_t^\alpha{f(g(t))}=f^\prime_g(g(t))D_t^\alpha{g},\mbox{(Chain rule)},\\
\label{eq1:6}
&D_t^\alpha{f}=t^{1-\alpha}f^\prime(t),f\,\mbox{is differentiable}.
\end{eqnarray}
Therefore, it is an acceptable choice to use CFD to solve the
conformable fractional partial differential equations (CFPDEs).
As a matter of fact,the CFD (\ref{eq1:1}) with its properties
(\ref{eq1:4}),(\ref{eq1:5}) and (\ref{eq1:6}) were successfully
used to solve some CFPDEs by means of the Riccati equation
expansion~\cite{Chen,Liu},the sine--cosine method~\cite{Darvishi},
the $\exp\left(\phi(\varepsilon)\right)$--expansion~\cite{Hosseini},
the modified Kudryashov method~\cite{Hosseini,Korkmaz} and the
Jacobi elliptic function expansion method\cite{Tasbozan}, etc.
However,we find that all these existing direct methods cannot
be used to find two types of exact traveling wave solutions of
CFPDEs when their original equations have the non--integrable
background.The first type solution,namely the exponential--Weierstrass
type solution,is expressed by the product of exponential
function and Weierstrass elliptic function. The second type
solution,namely the exponential--Jacobi type solution, is
expressed by the product of exponential function and Jacobi
elliptic function.Motivated by this problem,in this paper,
we shall propose a direct algebraic method for constructing
these two special types of solutions for CFPDEs.
It is seen that the method introduced in this paper is effective
to find the exponential--Weierstrass type,the exponential--Jacobi
type and the solitary type solutions of the CFPDEs.\par
This paper is organized as follows.In Sec.~\ref{sec:2},
by introducing two new second order auxiliary nonlinear ordinary
differential equations and their solutions we shall propose a new
algebraic method for solving CFPDEs.
In Sec.~\ref{sec:3},we shall use our suggested method to find
the exact traveling wave solutions of some CFPDEs. Finally,
Sec.~\ref{sec:4} offers a discussion.

\section{Description of the method}\label{sec:2}

Now let us simplify describe our direct algebraic method for finding
traveling wave solutions for CFPDEs. Suppose that a CFPDE is given by
\begin{equation}
\label{eq2:1}
P(u,D_t^\alpha{u},D_x^\beta{u},D_t^{2\alpha}u,{D_t^\alpha}D_x^\beta{u},
   D_x^{2\beta}u,\cdots)=0,
\end{equation}
where $P$ is a polynomial of its unknowns and $0<\alpha,\beta\leq{1}$.
The key steps of our method are outlined as follows.\par
{\em Step 1:}\,Making the wave transformation
\begin{equation}
\label{eq2:2}
u(x,t)=u(\xi),\xi={\frac {x^\beta}{\beta}}+\omega{\frac {t^\alpha}{\alpha}},
\end{equation}
and using the properties of CFD we can convert Eq.(\ref{eq2:1})
into the following ODE
\begin{equation}
\label{eq2:3}
Q(u,{\frac {du}{d\xi}},{\frac {d^2u}{d\xi^2}},\cdots)=0.
\end{equation}
\par
{\em Step 2:}\,(a)\,To seek the exponential--Weierstrass type
and solitary type exact traveling wave solutions,we
assume that Eq. (\ref{eq2:3}) has the solution of the form
\begin{equation}
\label{eq2:4}
u(\xi)=a_0+a_1F(\xi),
\end{equation}
where $a_0,a_1$ are constants, $F(\xi)$ satisfies the
following second order auxiliary ODE\cite{Siren}
\begin{equation}
\label{eq2:5}
F^{\prime\prime}(\xi)=bF^2(\xi)-6a^2F(\xi)+5aF^{\prime}(\xi).
\end{equation}
This equations admits the following solutions
\begin{equation}
\label{eq2:6}
F(\xi)=\left\{\begin{aligned}
&{\frac 6b}e^{2a\xi}\wp\left({\frac 1a}e^{a\xi}+c_1,0,g_3\right),\\
&{\frac {3a^2}{2b}}\left[1+\tanh\left({\frac a2}\xi\right)\right]^2,\\
&{\frac {3a^2}{2b}}\left[1+\coth\left({\frac a2}\xi\right)\right]^2,
\end{aligned}\right.
\end{equation}
(b)\,To find the exponential--Jacobi type and the solitary type
exact traveling wave solutions,we can take the solution of
Eq. (\ref{eq2:3}) of the form
\begin{equation}
\label{eq2:7}
u(\xi)=a_1F(\xi),
\end{equation}
where $a_1$ is an undermined constant, $F(\xi)$ satisfies the following
second order auxiliary ODE\cite{Siren}
\begin{equation}
\label{eq2:8}
F^{\prime\prime}(\xi)=cF^3(\xi)-2a^2F(\xi)-3aF^{\prime}(\xi).
\end{equation}
This equation has the following solutions
\begin{equation}
\label{eq2:9}
F(\xi)=\left\{\begin{array}{l}
{\varepsilon}ae^{-a\xi}\ds\left(e^{-a\xi}+c_2,{\frac {\sqrt{2}}2}\right),c=2,\\
{\varepsilon}ae^{-a\xi}\nc\left(\sqrt{2}e^{-a\xi}+c_2,{\frac {\sqrt{2}}2}\right),c=2,\\
{\frac {{\varepsilon}a}{2}}\left[1-\tanh\left({\frac a2}\xi\right)\right],c=2,\\
{\frac {{\varepsilon}a}{2}}\left[1-\coth\left({\frac a2}\xi\right)\right],c=2,\\
{\varepsilon}ae^{-a\xi}\cn\left(\sqrt{2}e^{-a\xi}+c_3,{\frac {\sqrt{2}}2}\right),c=-2,\\
{\frac {\sqrt{2}}2}{\varepsilon}ae^{-a\xi}\sd\left(\sqrt{2}e^{-a\xi}+c_3,{\frac {\sqrt{2}}2}\right),c=-2,
\end{array}\right.
\end{equation}
where $\wp$ expresses the Weierstras elliptic function,
$\ds,\nc,\cn,\sd$ are the Jacobi eliptic functions,$\varepsilon=\pm{1}$,
and $a,b,c_1,c_2,c_3,g_3$ are constants.\par
{\em Step 3:}\,Substituting (\ref{eq2:4}) with  (\ref{eq2:5}),
and (\ref{eq2:7}) with (\ref{eq2:8}) separately into Eq. (\ref{eq2:3})
and setting the coefficients of like powers of $F^i(F^{\prime})^j$ to zero,
we get a set of algebraic equations for unknowns $a_0,a_1,a,b,\omega$
,and for unknowns $a_1,a,\omega$,respectively. The system of algebraic
equations is solved by using a computer algebraic system,then the
values of these unknowns can be obtained.\par
{\em Step 4:}\,The exact traveling wave solutions of the
Eq.(\ref{eq2:1}) can be obtained by putting the values of
unknowns obtained in {\em Step 3} with (\ref{eq2:6}) into
(\ref{eq2:4}),and (\ref{eq2:9}) into (\ref{eq2:7}),respectively.

\section{Applications of the method}\label{sec:3}

Now we consider some illustrative examples to show the effectiveness
of our proposed method.
\begin{example}
\label{Ex1}
The space--time fractional KdV--Burgers equation\cite{El--Shewy}
\begin{align}
\label{eq3:1}
D_t^\alpha{u}+{\lambda}uD_x^\beta{u}+{\mu}{D_x^{2\beta}u}+\nu{D_x^{3\beta}u}=0,
\end{align}
in which $\lambda,\mu,\nu$ are constants.
\end{example}
Taking (\ref{eq2:2}) into (\ref{eq3:1}) we obtain the following ODE
\begin{equation}
\label{eq3:2}
{\omega}u^\prime+{\lambda}uu^\prime+{\mu}u^{\prime\prime}+{\nu}u^{\prime\prime\prime}=0.
\end{equation}
Substituting (\ref{eq2:4}) with (\ref{eq2:5}) into (\ref{eq3:2}) and setting
the coefficients of $F,F^2$,
$F^\prime,FF^\prime$ to be zero,we obtain a set of algebraic equations
\begin{eqnarray*}
\left\{\begin{array}{l}
2b{\nu}a_1+{\lambda}a_1^2=0,\\
-30a^3{\nu}a_1-6a^2{\mu}a_1=0,\\
5ab{\nu}a_1+b{\mu}a_1=0,\\
19a^2{\nu}a_1+5a{\mu}a_1+{\lambda}a_0a_1+{\omega}a_1=0.
\end{array}\right.
\end{eqnarray*}
It solves that
\begin{equation}
\label{eq3:3}
a=-{\frac {\mu}{5\nu}},
a_0={\frac {6\mu^2-25\omega\nu}{25\lambda\nu}},
a_1=-{\frac {2b\nu}{\lambda}}.
\end{equation}
Inserting (\ref{eq2:6}) with (\ref{eq3:3}) into (\ref{eq2:4}) we get
the exponential--Weierstrass type and solitary type exact traveling
wave solutions of Eq.(\ref{eq3:1}) as following
\begin{equation*}
\begin{array}{l}
u_1(x,t)=-{\frac {12}{\lambda}}e^{-\frac {2\mu}{5\nu}\left({\frac {x^\beta}{\beta}}+\omega{\frac {t^\alpha}{\alpha}}\right)}
\wp\left({\frac {5\nu}{\mu}}e^{-\frac {\mu}{5\nu}\left({\frac {x^\beta}{\beta}}
  +\omega{\frac {t^\alpha}{\alpha}}\right)}+c_1,0,g_3\right)
  +{\frac {6\mu^2-25\omega\nu}{25\lambda\nu}},\\
u_2(x,t)=-{\frac {3\mu^2}{25\lambda\nu}}\left[1-\tanh{\frac {\mu}{10\nu}}\left({\frac {x^\beta}{\beta}}
  +\omega{\frac {t^\alpha}{\alpha}}\right)\right]^2+{\frac {6\mu^2-25\omega\nu}{25\lambda\nu}},\\
u_3(x,t)=-{\frac {3\mu^2}{25\lambda\nu}}\left[1-\coth{\frac {\mu}{10\nu}}\left({\frac {x^\beta}{\beta}}
  +\omega{\frac {t^\alpha}{\alpha}}\right)\right]^2+{\frac {6\mu^2-25\omega\nu}{25\lambda\nu}},
  \end{array}
\end{equation*}
where $c_1,g_3$ are arbitrary constants.
\begin{example}\label{Ex2}
The time fractional Fisher equation\cite{Sungu}
\begin{equation}
 \label{eq3:4}
 D_t^\alpha{u}=u_{xx}+6u\left(1-u\right).
 \end{equation}
\end{example}
Substituting (\ref{eq2:2}) with $\beta=1$ into (\ref{eq3:4}),
we can convert the Eq. (\ref{eq3:4}) into the following ODE
\begin{equation}
\label{eq3:5}
{\omega}u^\prime=u^{\prime\prime}+6u\left(1-u\right),
\end{equation}
Substituting (\ref{eq2:4}) and (\ref{eq2:5}) into (\ref{eq3:5}) and
setting the coefficients of $F^\prime,F^j~(j=0,1,2)$ to zero,we have
\begin{eqnarray*}
\left\{\begin{array}{l}
6a_0^2-6a_0=0,\\
6a_1^2-a_1b=0,\\
-5aa_1+a_1\omega=0,\\
6a^2a_1+12a_0a_1-6a_1=0.
\end{array}\right.
\end{eqnarray*}
This algebraic equations is solved that
\begin{eqnarray}
\label{eq3:6}
&a_0=0,a_1={\frac b6},a=1,\omega=5,\\
\label{eq3:7}
&a_0=0,a_1={\frac b6},a=-1,\omega=-5.
\end{eqnarray}
Now the exponential--Weierstrass type and the solitary type
exact traveling wave solutions for Eq. (\ref{eq3:4}) can be
obtained by taking (\ref{eq2:6}) with (\ref{eq3:6}) and
(\ref{eq3:7}) into (\ref{eq2:4}),respectively,they are
\begin{equation*}
\begin{array}{l}
u_1(x,t)=e^{2x+{\frac {10t^\alpha}{\alpha}}}\wp\left(e^{x+{\frac {5t^\alpha}{\alpha}}}+c_1,0,g_3\right),\\
u_2(x,t)={\frac 14}\left[1+\tanh{\frac 12}\left(x
+{\frac {5t^\alpha}{\alpha}}\right)\right]^2,\\
u_3(x,t)={\frac 14}\left[1+\coth{\frac 12}\left(x
+{\frac {5t^\alpha}{\alpha}}\right)\right]^2,\\
u_4(x,t)=e^{-2x+{\frac {10t^\alpha}{\alpha}}}\wp\left(e^{-x+{\frac {5t^\alpha}{\alpha}}}+c_1,0,g_3\right),\\
u_5(x,t)={\frac 14}\left[1-\tanh{\frac 12}\left(x
-{\frac {5t^\alpha}{\alpha}}\right)\right]^2,\\
u_6(x,t)={\frac 14}\left[1-\coth{\frac 12}\left(x
-{\frac {5t^\alpha}{\alpha}}\right)\right]^2,
\end{array}
\end{equation*}
where $c_1,g_3$ are arbitrary constants.
 \begin{example}\label{Ex3}
The time fractional RLW--Burgers equation\cite{Korkmaz}
 \begin{equation}
 \label{eq3:8}
 D_t^\alpha{u}+pu_x+quu_x+ru_{xx}+su_{xxt}=0.
 \end{equation}
 \end{example}
Inserting (\ref{eq2:2}) with $\beta=1$ into (\ref{eq3:8}) leads the
following ODE
\begin{equation}
\label{eq3:9}
\left(\omega+p\right)u^\prime+quu^\prime+ru^{\prime\prime}+{\omega}su^{\prime\prime\prime}=0.
\end{equation}
Taking (\ref{eq2:4}) with (\ref{eq2:5}) into (\ref{eq3:9}) and setting
the coefficients of $F,F^2,F^\prime,FF^\prime$ to be zero,then we obtain
\begin{equation*}
\left\{\begin{array}{l}
-30a^3{\omega}sa_1-6a^2ra_1=0,\\
2b{\omega}sa_1+qa_1^2=0,\\
5ab{\omega}sa_1+bra_1=0,\\
19a^2{\omega}sa_1+5ara_1+qa_0a_1+\left(\omega+p\right)a_1=0.
\end{array}\right.
\end{equation*}
It solves that
\begin{equation}
\label{eq3:10}
a=-{\frac r{5s\omega}},
a_0=-{\frac {25s\omega^2+25ps\omega-6r^2}{25qs\omega}},
a_1=-{\frac {2bs\omega}q}.
\end{equation}
The exponential--Weierstrass type and the solitary type exact solutions
of the Eq.(\ref{eq3:8}) are obtained by substituting (\ref{eq2:6}) with
(\ref{eq3:10}) into (\ref{eq2:4}),they are the following
\begin{equation*}
\begin{array}{l}
u_1(x,t)=-{\frac {12}q}e^{-\frac {2r}{5s\omega}\left(x+{\frac {{\omega}t^\alpha}{\alpha}}\right)}
\wp\left({\frac {5s\omega}r}e^{-\frac {r}{5s\omega}\left(x+{\frac {{\omega}t^\alpha}{\alpha}}\right)}+c_1,0,g_3\right)\\
\qquad-{\frac {25s\omega^2+25ps\omega-6r^2}{25qs\omega}},\\
u_2(x,t)=-{\frac {3r^2}{25qs\omega}}\left[1-\tanh{\frac {r}{10s\omega}}\left(x+{\frac {\omega{t^\alpha}}{\alpha}}\right)\right]^2
  -{\frac {25s\omega^2+25ps\omega-6r^2}{25qs\omega}},\\
u_3(x,t)=-{\frac {3r^2}{25qs\omega}}\left[1-\coth{\frac {r}{10s\omega}}\left(x+{\frac {\omega{t^\alpha}}{\alpha}}\right)\right]^2
  -{\frac {25s\omega^2+25ps\omega-6r^2}{25qs\omega}},
\end{array}
\end{equation*}
where $c_1,g_3$ are free parameters.
\begin{example}\label{Ex4}
The time fractional Cahn--Allen equation\cite{Khater}
\begin{equation}
 \label{eq3:11}
 D_t^\alpha{u}-u_{xx}-u+u^3=0.
\end{equation}
\end{example}
Substituting (\ref{eq2:2}) with $\beta=1$ into (\ref{eq3:11})
we get the following ODE
\begin{equation}
 \label{eq3:12}
 {\omega}u^\prime-u^{\prime\prime}-u+u^3=0.
\end{equation}
 Substituting (\ref{eq2:7}) with (\ref{eq2:8}) into (\ref{eq3:12})
 and setting the coefficients of $F,F^3,F^\prime$ to zero gives the
 following ODE
\begin{equation*}
 \left\{\begin{array}{l}
 2a^2a_1-a_1=0,\\
 a_1^3-ca_1=0,\\
 3aa_1+{\omega}a_1=0.
 \end{array}\right.
\end{equation*}
Solving this algebraic equation we obtain that
\begin{equation}
 \label{eq3:13}
 a=\pm{\frac {\sqrt{2}}2},a_1=\pm\sqrt{c},\omega=\mp{\frac {3\sqrt{2}}2}.
\end{equation}
The exponential--Jacobi type and the solitary type exact solutions
of Eq. (\ref{eq3:11}) are obtained by taking (\ref{eq3:13}) and
(\ref{eq2:9}) with $c=2$ into (\ref{eq2:7}),they are now given by
\begin{eqnarray*}
 \begin{array}{l}
 u_1(x,t)={\varepsilon}e^{\mp{\frac {\sqrt{2}}2}\left(x\mp{\frac {3\sqrt{2}}{2\alpha}t^\alpha}\right)}
   \ds\left(e^{\mp{\frac {\sqrt{2}}2}\left(x\mp{\frac {3\sqrt{2}}{2\alpha}t^\alpha}\right)}+c_2,{\frac {\sqrt{2}}2}\right),\\
 u_2(x,t)={\varepsilon}e^{\mp{\frac {\sqrt{2}}2}\left(x\mp{\frac {3\sqrt{2}}{2\alpha}t^\alpha}\right)}
 \nc\left(\sqrt{2}e^{\mp{\frac {\sqrt{2}}2}\left(x\mp{\frac {3\sqrt{2}}{2\alpha}t^\alpha}\right)}+c_2,{\frac {\sqrt{2}}2}\right),\\
 u_3(x,t)={\frac {\varepsilon}2}\left[1\mp\tanh{\frac {\sqrt{2}}{4}}\left(x\mp{\frac {3\sqrt{2}t^\alpha}{2\alpha}}\right)\right],\\
 u_4(x,t)={\frac {\varepsilon}2}\left[1\mp\coth{\frac {\sqrt{2}}{4}}\left(x\mp{\frac {3\sqrt{2}t^\alpha}{2\alpha}}\right)\right],
 \end{array}
\end{eqnarray*}
where $c_2$ is a free parameter.
\begin{example}\label{Ex5}
The space--time fractional mKdV--Burgers equation~\cite{Liu}
\begin{equation}
 \label{eq3:14}
 D_t^\alpha{u}+{\lambda}u^2D_x^\alpha{u}+rD_x^{2\alpha}u+sD_x^{3\alpha}u=0.
 \end{equation}
\end{example}
 Inserting (\ref{eq2:2}) with $\beta=\alpha$ into (\ref{eq3:14}) leads
 the following ODE
\begin{equation}
 \label{eq3:15}
 {\omega}u^\prime+{\lambda}u^2u^\prime+ru^{\prime\prime}+su^{\prime\prime\prime}=0.
\end{equation}
Taking (\ref{eq2:7}) with (\ref{eq2:8}) into (\ref{eq3:15}) and setting
the coefficients of $F,F^2,F^\prime,F^2F^\prime$ to zero,we obtain
a set of algebraic equations
\begin{equation*}
 \left\{\begin{array}{l}
 6a^3sa_1-2a^2ra_1=0,\\
 {\lambda}a_1^3+3csa_1=0,\\
 -3acsa_1+cra_1=0,\\
 7a^2sa_1-3ara_1+{\omega}a_1=0.
 \end{array}\right.
\end{equation*}
 Its solution is found to be
\begin{equation}
 \label{eq3:16}
 a={\frac r{3s}},a_1=\pm\sqrt{-\frac {3cs}{\lambda}},\omega={\frac {2r^2}{9s}}.
\end{equation}
The exponential--Jacobi type and the solitary type exact traveling
wave solutions of Eq.(\ref{eq3:14}) obtained by taking
(\ref{eq3:16}) and (\ref{eq2:9}) with $c=2$ into (\ref{eq2:7}) are
\begin{eqnarray*}
 \begin{array}{l}
 u_1(x,t)={\frac {r\varepsilon}{3s}}\sqrt{-\frac {6s}{\lambda}}e^{-{\frac {r}{3s}}\left({\frac {x^\alpha}{\alpha}}+{\frac {2r^2t^\alpha}{9s\alpha}}\right)}
 \ds\left(e^{-{\frac {r}{3s}}\left({\frac {x^\alpha}{\alpha}}+{\frac {2r^2t^\alpha}{9s\alpha}}\right)}+c_2,{\frac {\sqrt{2}}2}\right),s\lambda<0,\\
 u_2(x,t)={\frac {r\varepsilon}{3s}}\sqrt{-\frac {6s}{\lambda}}e^{-{\frac {r}{3s}}\left({\frac {x^\alpha}{\alpha}}+{\frac {2r^2t^\alpha}{9s\alpha}}\right)}
 \nc\left(\sqrt{2}e^{-{\frac {r}{3s}}\left({\frac {x^\alpha}{\alpha}}+{\frac {2r^2t^\alpha}{9s\alpha}}\right)}+c_2,{\frac {\sqrt{2}}2}\right),s\lambda<0,\\
 u_3(x,t)=-{\frac {r\varepsilon}{6s}}\sqrt{-\frac {6s}{\lambda}}
   \left[1-\tanh{\frac r{6s}}\left({\frac {x^\alpha}{\alpha}}+{\frac {2r^2t^\alpha}{9s\alpha}}\right)\right],s\lambda<0,\\
 u_4(x,t)=-{\frac {r\varepsilon}{6s}}\sqrt{-\frac {6s}{\lambda}}
   \left[1-\coth{\frac r{6s}}\left({\frac {x^\alpha}{\alpha}}+{\frac {2r^2t^\alpha}{9s\alpha}}\right)\right],s\lambda<0,
   \end{array}
\end{eqnarray*}
where $c_2$ is an arbitrary constant.\par
By inserting (\ref{eq3:16}) and (\ref{eq2:9}) with $c=-2$
into (\ref{eq2:7}) we get the following exponential--Jacobi type and
solitary type traveling wave solutions of Eq.(\ref{eq3:11})
\begin{eqnarray*}
 u_5(x,t)={\frac {r\varepsilon}{3s}}\sqrt{\frac {6s}{\lambda}}e^{-{\frac {r}{3s}}\left({\frac {x^\alpha}{\alpha}}+{\frac {2r^2t^\alpha}{9s\alpha}}\right)}
 \cn\left(\sqrt{2}e^{-{\frac {r}{3s}}\left({\frac {x^\alpha}{\alpha}}+{\frac {2r^2t^\alpha}{9s\alpha}}\right)}+c_2,{\frac {\sqrt{2}}2}\right),s\lambda>0,\\
 u_6(x,t)={\frac {r\varepsilon}{3s}}\sqrt{\frac {3s}{\lambda}}e^{-{\frac {r}{3s}}\left({\frac {x^\alpha}{\alpha}}+{\frac {2r^2t^\alpha}{9s\alpha}}\right)}
 \sd\left(\sqrt{2}e^{-{\frac {r}{3s}}\left({\frac {x^\alpha}{\alpha}}+{\frac {2r^2t^\alpha}{9s\alpha}}\right)}+c_2,{\frac {\sqrt{2}}2}\right),s\lambda<0,
\end{eqnarray*}
in which $c_2$ is a free parameter.
\begin{example}\label{Ex6}
The space--time fractional telegraph equation\cite{Liu}
\begin{equation}
 \label{eq3:17}
 D_t^{2\alpha}u-D_x^{2\alpha}u+D_t^\alpha{u}+{\mu}u+{\nu}u^3=0.
\end{equation}
\end{example}
Taking (\ref{eq2:2}) with $\beta=\alpha$ into (\ref{eq3:17})
we obtain the following ODE
\begin{equation}
 \label{eq3:18}
 \left(\omega^2-1\right)u^{\prime\prime}+{\omega}u^\prime+{\mu}u+{\nu}u^3=0.
\end{equation}
Substituting (\ref{eq2:7}) with (\ref{eq2:8}) into (\ref{eq3:18}) and
setting the coefficients of $F,F^3,F^\prime$ to zero we obtain
\begin{equation*}
 \left\{\begin{array}{l}
 -3a\omega^2a_1+3aa_1+{\omega}a_1=0,\\
 -2a^2\omega^2a_1+2a^2a_1+{\mu}a_1=0,\\
 {\nu}a_1^3+c\omega^2a_1-ca_1=0.
 \end{array}\right.
\end{equation*}
Solutions of this algebraic equations are found to be
\begin{eqnarray}
 \label{eq3:19}
 a={\frac {9\mu-2}2}\sqrt{\frac {\mu}{9\mu-2}},
 a_1=\pm\sqrt{-\frac {2c}{\nu\left(9\mu-2\right)}},
 \omega=3\sqrt{\frac {\mu}{9\mu-2}},\\
 \label{eq3:20}
 a=-{\frac {9\mu-2}2}\sqrt{\frac {\mu}{9\mu-2}},
 a_1=\pm\sqrt{-\frac {2c}{\nu\left(9\mu-2\right)}},
 \omega=-3\sqrt{\frac {\mu}{9\mu-2}},
\end{eqnarray}
When $c=2$ and $c=-2$,by substituting (\ref{eq3:19}) and (\ref{eq2:9})
into (\ref{eq2:7}),respectively,we get the exponential--Jacobi type
and the solitary type traveling wave solutions of Eq.(\ref{eq3:17})
as following
\begin{equation*}
\begin{array}{l}
u_1(x,t)=\varepsilon\sqrt{-\frac {\mu}{\nu}}e^{-\eta^{+}}
 \ds\left(e^{-\eta^{+}}+c_2,{\frac {\sqrt{2}}{2}}\right),
 \mu<0,\nu>0,\mbox{or}\;\mu>{\frac 29},\nu<0,\\
u_2(x,t)=\varepsilon\sqrt{-\frac {\mu}{\nu}}e^{-\eta^{+}}
 \nc\left(\sqrt{2}e^{-\eta^{+}}+c_2,{\frac {\sqrt{2}}{2}}\right),\mu<0,\nu>0,
 \mbox{or}\;\mu>{\frac 29},\nu<0,\\
u_3(x,t)={\frac {\varepsilon}2}\sqrt{-\frac {\mu}{\nu}}
  \left[1-\tanh\eta^{+}\right],\mu<0,\nu>0,\mbox{or}\;\mu>{\frac 29},\nu<0,\\
u_4(x,t)={\frac {\varepsilon}2}\sqrt{-\frac {\mu}{\nu}}
  \left[1-\coth\eta^{+}\right],\mu<0,\nu>0,\mbox{or}\;\mu>{\frac 29},\nu<0,\\
u_5(x,t)=\varepsilon\sqrt{\frac {\mu}{\nu}}e^{-\eta^{+}}
 \cn\left(\sqrt{2}e^{-\eta^{+}}+c_2,{\frac {\sqrt{2}}{2}}\right),\mu<0,\nu<0,\mbox{or}\;\mu>{\frac 29},\nu>0,\\
u_6(x,t)={\frac {\varepsilon}2}\sqrt{\frac {2\mu}{\nu}}e^{-\eta^{+}}
 \sd\left(\sqrt{2}e^{-\eta^{+}}+c_2,{\frac {\sqrt{2}}{2}}\right),\mu<0,\nu<0,\mbox{or}\;\mu>{\frac 29},\nu>0,
\end{array}
\end{equation*}
where
\begin{equation*}
\eta^{+}={\frac {9\mu-2}2}\sqrt{\frac {\mu}{9\mu-2}}\left({\frac {x^\alpha}{\alpha}}+{\frac {3}{\alpha}}\sqrt{{\frac {\mu}{9\mu-2}}}t^\alpha\right).
\end{equation*}
When $c=2$ and $c=-2$,by substituting (\ref{eq3:20}) and (\ref{eq2:9})
into (\ref{eq2:7}),respectively,the exponential--Jacobi type
and solitary type traveling wave solutions of Eq.(\ref{eq3:17}) are
obtained as
\begin{equation*}
\begin{array}{l}
u_7(x,t)=\varepsilon\sqrt{-\frac {\mu}{\nu}}e^{\eta^{-}}
 \ds\left(e^{\eta^{-}}+c_2,{\frac {\sqrt{2}}{2}}\right),\mu<0,\nu>0\;\mbox{or}\;\mu>{\frac 29},\nu<0,\\
u_8(x,t)=\varepsilon\sqrt{-\frac {\mu}{\nu}}e^{\eta^{-}}
 \nc\left(\sqrt{2}e^{\eta^{-}}+c_2,{\frac {\sqrt{2}}{2}}\right),\mu<0,\nu>0\;\mbox{or}\;\mu>{\frac 29},\nu<0,\\
u_9(x,t)={\frac {\varepsilon}2}\sqrt{-\frac {\mu}{\nu}}
  \left[1+\tanh\eta^{-}\right],\mu<0,\nu>0,\mbox{or}\;\mu>{\frac 29},\nu<0,\\
u_{10}(x,t)={\frac {\varepsilon}2}\sqrt{-\frac {\mu}{\nu}}
  \left[1+\coth\eta^{-}\right],\mu<0,\nu>0,\mbox{or}\;\mu>{\frac 29},\nu<0,\\
u_{11}(x,t)=\varepsilon\sqrt{\frac {\mu}{\nu}}e^{\eta^{-}}
 \cn\left(\sqrt{2}e^{\eta^{-}}+c_2,{\frac {\sqrt{2}}{2}}\right),\mu<0,\nu<0\;\mbox{or}\;\mu>{\frac 29},\nu>0,\\
u_{12}(x,t)={\frac {\varepsilon}2}\sqrt{\frac {2\mu}{\nu}}e^{\eta^{-}}
 \sd\left(\sqrt{2}e^{\eta^{-}}+c_2,{\frac {\sqrt{2}}{2}}\right),\mu<0,\nu<0\;\mbox{or}\;\mu>{\frac 29},\nu>0,
 \end{array}
 \end{equation*}
where
\begin{equation*}
\eta^{-}={\frac {9\mu-2}2}\sqrt{\frac {\mu}{9\mu-2}}\left({\frac {x^\alpha}{\alpha}}-{\frac {3}{\alpha}}\sqrt{{\frac {\mu}{9\mu-2}}}t^\alpha\right).
\end{equation*}
\section{Discussion}
\label{sec:4}
The examples in Sec.~\ref{sec:3} showed that our method is very effective
to find the exponential--Weierstrass type,the exponential--Jacobi type and
solitary type exact traveling wave solutions for CFPDEs.
In our method, the auxiliary equations (\ref{eq2:5}) and (\ref{eq2:8})
are new and they have not been previously used in any direct methods.
Therefore,all our exponential--Weierstrass type and exponential--Jacobi
type traveling wave solutions are new and cannot be found by using the
existing direct algebraic methods.
It is also pointed out that the original integer order nonlinear
partial differential equations of the CFPDEs considered in our
examples are non--integrable. Therefore,our method is restricted
to solve those CFPDEs whose original equations are non--integrable.

\section*{Acknowledgements}

This work has been supported by the National Natural Science Foundation of
China under Grant No.11861050.


\end{document}